\documentclass{article}
\usepackage{latexsym}
\usepackage{amsfonts}
\usepackage{amssymb}
\usepackage{color}

\newtheorem{theorem}{Theorem}[section]
\newtheorem{lemma}[theorem]{Lemma}
\newtheorem{corollary}[theorem]{Corollary}
\newtheorem{proposition}[theorem]{Proposition}
\newtheorem{claim}[theorem]{Claim}
\newtheorem{question}[theorem]{Question}
\newtheorem{problem}[theorem]{Problem}

\newtheorem{definition}[theorem]{Definition}
\newtheorem{remark}[theorem]{Remark}

\begin{document}

\newcommand{\normal}{\vartriangleleft}
\newcommand{\nnormal}{\ntriangleleft}
\newcommand{\R}{\mathbb{R}}
\newcommand{\Z}{\mathbb{Z}}
\newcommand{\F}{\mathbb{F}}
\newcommand{\A}{\mathcal{A}}
\newcommand{\bbC}{\mathbb{C}}
\newcommand{\abs}[1]{\left| #1\right|}
\newcommand{\ord}[1]{\Delta \left( #1 \right)}
\newcommand{\leqs}{\leqslant}
\newcommand{\geqs}{\geqslant}
\newcommand{\heq}{\simeq}
\newcommand{\iso}{\simeq}
\newcommand{\maps}{\longrightarrow}
\newcommand{\meet}{\wedge}
\newcommand{\join}{\vee}
\newcommand{\homeo}{\cong}
\newcommand{\isom}{\cong}
\newcommand{\cross}{\times}
\newcommand{\mC}{\mathcal{C}}
\newcommand{\C}[1]{\mathcal{C}(#1)}
\newcommand{\M}[1]{\mathcal{M}(#1)}
\newcommand{\supp}[1]{\rm supp\left( #1 \right)}
\newcommand{\gen}[1]{\langle #1 \rangle}
\newcommand{\prj}{\F_p \cup \{\infty\}}
\newcommand{\gl}{GL_2 (\F_p)}
\newcommand{\spl}{SL_2 (\F_p)}
\newcommand{\psl}{PSL_2 (\F_p)}
\newcommand{\psls}{PSL_2 (\F_7)}
\newcommand{\injects}{\hookrightarrow}
\newcommand{\mF}{\mathcal{F}}

\title{Connectivity of the coset poset and the subgroup poset of a group}

\author{Daniel A. Ramras}

\date{}

\maketitle{}

\begin{abstract} We study the connectivity of the coset poset and the subgroup
poset of a group, focusing in particular on simple connectivity.
The coset poset was recently introduced by K. S. Brown in connection with
the probabilistic zeta function of a group.
We further Brown's study of the homotopy type of the coset poset, and
in particular generalize
his results on direct products and classify
direct products with simply connected coset posets.

The homotopy type of the subgroup poset $L(G)$ has been examined previously 
by Kratzer, Th\'{e}venaz, and 
Shareshian.  We generalize some results of Kratzer and Th\'{e}venaz, and determine 
$\pi_1 (L(G))$ in nearly all cases.
\end{abstract}  

\section{Introduction}$\label{intro}$

One may apply topological concepts to any poset (partially ordered set) $P$ by means
of the simplicial complex $\ord{P}$ consisting of all finite chains of $P$.  The topology
of posets arising from groups has been studied extensively; 
see for 
example~\cite{Brown-coset-poset,Kratzer-subgroups, Quillen-p-subgroups, Thevenaz-top-homology}.
The basic topological theory of posets is described in~\cite{Bjorner-top-methods} 
and in~\cite{Ramras-thesis}.
We will often use $P$ to denote both $\ord{P}$ and its geometric realization 
$\abs{\ord{P}} = \abs{P}$.

The coset poset $\C{G}$ of a finite group $G$ (the poset of all left cosets of all proper 
subgroups of $G$, ordered by inclusion) was introduced by Brown~\cite{Brown-coset-poset} 
in connection with the probabilistic zeta function $P(G, s)$.
(Note that the choice of left cosets over right cosets is irrelevant since each left coset
is a right coset ($xH = (xHx^{-1}) x$) and vice versa.)
Brown showed that $P(G, -1) = -\tilde{\chi} (\ord{\C{G}})$ and used this relationship to
prove certain divisibility results about $P(G,-1)$ (which is always an integer). 
In fact, Brown proved that for arbitrary integers $s$, $P(G,s)$ can be 
calculated from the coset poset.  The connection between $P(G, -1)$ and $\chi (\C{G})$ motivates
the study of the homotopy type of $\C{G}$.  Brown obtained a variety of results along these
lines, and in this paper we prove a number of new results about the homotopy-type of $\C{G}$. 
In particular, we study the connectivity of $\C{G}$ and further Brown's study of the homotopy
type of $\C{G}$ in terms of normal subgroups and quotients.  Brown
has asked~\cite[Question 3]{Brown-coset-poset}, ``For which finite groups $G$ 
is $\C{G}$ simply connected?''  One of the 
main goals of this paper is to study this question (for both finite and infinite groups).  
We now gather together our main results on the problem.

\begin{enumerate}
	\item If $G$ is not a 2-generator group, then $\C{G}$ is simply connected
		(Corollary~\ref{2-generator}).
	\item If $G = H\cross K$ ($H, K$ non-trivial), then $\pi_1 (\C{G})\neq 1$ if and only if
		both $H$ and $K$ are cyclic of prime-power order (Theorem~\ref{direct-products}).
	\item If $G = H\rtimes K$ ($H, K$ non-trivial) where $K$ is not finite cyclic 
		or $H$ is not $K$-cyclic, then $\C{G}$ is simply connected 
		(Proposition~\ref{pi-semi-direct}).
	\item If the subgroup poset $L(G)$ is disconnected, 
		or if $G$ has a maximal subgroup isomorphic to 
		$\Z/p^n$ ($p$ prime), then $\C{G}$ is not simply connected 
		(in fact, lower bounds on the rank of $H_1 (\C{G})$ are obtained in 
		Theorems~\ref{M-V1} and~\ref{M-V2}).
	\item If $G$ is isomorphic to $PSL_2 (F_5)\isom A_5$ or $\psls$, then $\C{G}$ is simply
		connected (Propositions~\ref{A_5} and~\ref{psls-simply-connected}).
\end{enumerate}
Further results appear in Corollary~\ref{2-gen-presentation} and
Propositions~\ref{pi-semi-direct},~\ref{S_n}, and~\ref{coset-quotient}.

Brown has given a complete description of the homotopy type of $\C{G}$ for any
finite solvable group $G$~\cite[Proposition 11]{Brown-coset-poset}, 
and we now take a moment to discuss this result.  
Recall that a \emph{chief series} in a group $G$ is a maximal chain in the lattice of normal 
subgroups of $G$. 

\begin{theorem}[Brown]$\label{solvable-coset-poset}$
Let $G$ be a finite solvable group, let 
$\{1\} = N_0 \normal \cdots \normal N_k = G$
be a chief series for $G$ and let $c_i$ be the number of complements of $N_i/N_{i-1}$ in 
$G/N_{i-1}$ ($1\leqs i\leqs k$).  Then $\C{G}$ is homotopy equivalent to a bouquet of 
spheres, each of dimension $d-1$ where $d$ is the number of indices $\, i$ 
such that $c_i\neq 0$.  The number of spheres is given by 
$\abs{ \prod_{i=1}^{k} 
		\left( c_i \frac{\abs{N_i}
                                }
                                {\abs{N_{i-1}
                                     }
                                }
                       - 1 
                \right)
    }$.
\end{theorem}
The number of spheres may be calculated by induction, using~\cite[Corollary 3]{Brown-coset-poset}.
We remark that when $G$ is solvable, the number $d$ of complemented factors in a chief series for 
$G$ is bounded below by $\pi (G)$, the number of distinct primes dividing $o(G)$.  
This follows from the fact that any solvable group $G$ may be built up by a series of extensions 
with kernel of prime-power order.  In this process, each time we add a new prime the 
kernel $P$ and quotient $Q$ of the extension are relatively prime and hence 
$H^2 (Q, P) = 0$~\cite[Chapter 3, Proposition 10.1]{Brown-cohomology} and the extension splits.
In light of the above theorem, this shows that if $G$ is a finite solvable group, 
$\C{G}$ is $k$-connected whenever $\pi (G)\geqs k+2$.

The paper is structured as follows.  In Section~\ref{connectivity}, we introduce
the notion of an \emph{atomized} poset.  The coset poset of any group and the subgroup poset
of any finite group (or any torsion group) are atomized, and the main theorem of this section gives 
conditions under which an atomized poset is $k$-connected.  

The main results on the simple connectivity of $\C{G}$ appear in Section~\ref{extensions}.
We study the coset poset of a semi-direct product $H\rtimes K$ and extend
Brown's analysis of direct products to this setting.  

In Section~\ref{m-v}, we study the homology of $\C{G}$ using Mayer-Vietoris sequences,
and obtain conditions under which the coset poset is not simply connected.  

The results described so far have little bearing on finite, non-abelian simple groups.
In Section~\ref{simple-groups}, we treat the first two non-abelian simple groups and
prove that both have simply connected coset posets. 

Many of the techniques in this paper apply to the subgroup poset $L(G)$ (the
poset of all proper, non-trivial subgroups of $G$, ordered by inclusion) as well as
to the coset poset, and in Section~\ref{L(G)} we discuss the connectivity of $L(G)$.
The subgroup poset of a finite group has been studied previously by Kratzer and 
Th\'{e}venaz~\cite{Kratzer-subgroups} and by Th\'{e}venaz~\cite{Thevenaz-top-homology}.  
Using the homotopy complementation
formula of Bj\"{o}rner and Walker~\cite{Bjorner-hom-comp}, we generalize some of the results 
in~\cite{Kratzer-subgroups}.  Our results on the connectivity of $L(G)$ are based mainly on
this complementation formula and on an analysis of the connectedness of $L(G)$.  Some
of the results could have been proven using methods from earlier sections,
but the chosen approach is more elegant and more powerful.  
After proving general results about the connectivity of $L(G)$, we restrict our attention to
finite groups and determine $\pi_1 (L(G))$ in nearly all cases.

In the final section of the paper, we discuss a relationship between the
homology of the coset poset and that of the subgroup poset of certain groups.  Kratzer and 
Th\'{e}venaz~\cite{Kratzer-subgroups} have proven an analogue of Theorem~\ref{solvable-coset-poset}
for the subgroup poset of a finite solvable group (Theorem~\ref{solvable-subgroup-poset}).  
The striking similarity between these two results is part of the motivation behind our
discussion, and in particular shows that the relationship holds for finite solvable groups.

Some remarks on notation and conventions: If $\Delta$ is a simplicial complex, we will denote its
$k$-skeleton by $\Delta^{\leqs\, k}$.  In a poset $P$, a chain $p_1<p_2<\cdots <p_n$ is said
to have length $n-1$.  A wedge sum of empty topological spaces, or a wedge sum over an empty
index set, is defined to be a point.  Finally, we consider the empty space to be path-connected
but not 0-connected.

The research described in this article began as an undergraduate thesis~\cite{Ramras-thesis},
and certain results are described there in more detail 
(especially the results in Section~\ref{simple-groups}).
I thank Kenneth S. Brown, who served as my thesis advisor, for all his help in connection with 
this article.  
In addition to introducing me to the coset poset, Professor Brown suggested many of the 
techniques in this paper and helped a great deal with the exposition.  I also thank
Maria Silvia Lucido for providing the answer to a question I asked about the subgroup poset.

\section{Connectivity of Atomized Posets}$\label{connectivity}$

In this section we introduce \emph{atomized} posets, a class of posets generalizing 
the coset poset of a group.  For finite groups (and in fact torsion groups) the subgroup poset
is also atomized.  The main result of this section
is Theorem~\ref{k-connected-atom}, which gives conditions under which an atomized poset
is $k$-connected.

\begin{definition}  Let $P$ be a poset.  We call $P$ \emph{atomized} if every element of $P$ lies
above some minimal element and every finite set of minimal elements with an upper bound has a join.

We call the minimal elements of $P$ \emph{atoms}, and denote the set of atoms of $P$
by $\A(P)$.  If $S\subset \A (P)$ is finite,
we say that $S$ \emph{generates} its join  
(or that $S$ generates $P$ if $P_{\geqs S}$ is empty), and we write $\gen{S}$
for the object generated by $S$ (so $\gen{S} = P$ is allowed).
\end{definition}

The proper part of any finite-length lattice is atomized, but the converse
is not true.  (Consider, for example, the five-element poset $\{\hat{0}<a,b<c,d\}$.  This poset
has a unique minimal element and is thus atomized, but its maximal elements do not have a meet.)
Note that the subgroup poset of $\Z$ has no minimal elements and is not atomized, 
and in fact if a group $G$ has an element of infinite order then $L(G)$ is not atomized.

In the coset poset of any group or the subgroup poset of any torsion group,
the definition of generation coincides with the standard group theoretic definitions, where the
coset generated by elements $x_1, x_2, \ldots \in G$ 
is $$x_1\gen{x_1^{-1} x_2, x_1^{-1} x_3, \ldots }.$$

The following lemma shows 
that, up to homotopy, we can replace any atomized poset $P$ with a smaller simplicial complex, 
$\M{P}$.  
This complex has many fewer vertices, although it has much higher dimension, and will play
an important role in our analysis of the coset poset.

\begin{lemma} $\label{min-cover}$
Let $P$ be an atomized poset and let $\M{P}$ denote the simplicial complex with vertex set
$\A (P)$ and with a simplex for each finite set $S\subset \A (P)$ with $\gen{S} \neq P$.  
Then $\ord{P} \heq \M{P}$.
\end{lemma}
{\bf Proof.}  We will apply the Nerve Theorem~\cite[10.6]{Bjorner-top-methods}.
Consider the collection of all cones $P_{\geqs x}$ with $x \in \A (P)$.  Since
$P$ is atomized, $\bigcup_{x\in \A (P)} \ord{P_{\geqs x}} = \ord{P}$.
If $S\subset \A (P)$ is finite and $\gen{S}\neq P$, then
$\cap_{s\in S} P_{\geqs s} = P_{\geqs \gen{S}} \heq \{\ast\}.$  So each finite intersection 
is either empty or contractible, and the Nerve Theorem tells us that $\ord{P}$ is 
homotopy equivalent to the nerve of this cover, which is exactly $\M{P}$.
$\hfill \Box$
\vspace{.1in}

We call $\M{P}$ the minimal cover of $P$.
When $P = \C{G}$ for some group $G$, we denote $\M{\C{G}}$ by $\M{G}$.  This complex
has vertex set $G$ and a simplex for each finite subset of $G$ contained in a proper coset.

\begin{theorem}$\label{k-connected-atom}$
Let $P$ be an atomized poset such that no $k$ atoms generate $P$.  Then
$P$ is $(k-2)$-connected.  In particular, for any group $G$
\begin{enumerate}
	\item if $G$ is not $k$-generated, then $\C{G}$ is $(k-1)$-connected;
	\item if $G$ is a torsion group in which any $k$ elements of prime order generate
		a proper subgroup, then $L(G)$ is $(k-2)$-connected.
\end{enumerate}
\end{theorem}
{\bf Proof.}  By the lemma, it suffices to check that $\M{P}$ is $(k-2)$-connected.
Since no $k$ atoms generate $P$, any $k$ atoms form a simplex in $\M{P}$.  So
$\M{P}^{\leqs\, k-1}$ is just the $(k-1)$-skeleton of the full simplex on the vertex set $\A (P)$,  
and thus $\M{P}$ is $(k-2)$-connected.  (Note
that $\A (P)$ may be infinite; in this case the ``full simplex" on the set $\A (P)$ is the
simplicial complex whose simplices are all the finite subsets of $\A (P)$.)  
$\hfill \Box$

\begin{remark} Brown has asked~\cite[Question 4]{Brown-coset-poset} whether there exist
finite groups $G$ with $\C{G}$ contractible.  The above theorem shows that if $G$ is
a (necessarily infinite) group which is not finitely generated, then $\C{G}$ is contractible.
\end{remark}

Theorem~\ref{k-connected-atom} specializes to the following result giving conditions
under which $\C{G}$ (or any atomized poset) is simply connected.

\begin{corollary}$\label{2-generator}$
If $G$ is not a 2-generator group, then
$\C{G}$ is simply connected.  More generally, if $P$ is an atomized poset and
no three atoms generate $P$, then $\abs{P}$ is simply connected.
\end{corollary}

This theorem does not characterize finite groups with
simply connected coset posets.  In fact, the first non-abelian simple group, $A_5$, affords
a counterexample (Proposition~\ref{A_5}).

We now give an alternate proof of Corollary~\ref{2-generator}.  This proof will rely solely 
on the combinatorial structure of $\ord{P}$ itself (as opposed to that of $\M{P}$). 
First we need a lemma giving a simple representative
for each homotopy class of loops in $\abs{P}$.  Recall that an edge cycle $C$ in a 
graph $\Gamma$ is determined by listing (in order) the vertices through which $C$ passes.

\begin{lemma}$\label{cyclic-cycle}$
Let $P$ be an atomized poset in which no two elements generate.  
Then every loop in $\abs{P}$ is homotopic to an edge 
cycle in $\ord{P}^{\leqs\, 1}$ of the form 
$$(a_1, \gen{a_1, a_2}, a_2,\gen{a_2, a_3},\ldots, a_n, \gen{a_n, a_1}, a_1),$$ 
where the $a_i$ are atoms. 
\end{lemma}
{\bf Proof.}  By the Simplicial Approximation Theorem, every loop in $\abs{P}$ is homotopic
to an edge cycle.  Let $C = (p_1, p_2, \ldots , p_n = p_1)$ be an edge cycle in 
$\ord{P}^{\leqs\, 1}$.  
For each $i$ either $p_i < p_{i+1}$ or $p_i > p_{i+1}$, and it is easy to check that $C$ is 
homotopic to a cycle $C' = (p'_1,\ldots, p'_m)$
in which these inequalities alternate direction (simply remove the points at which the 
inequalities do not alternate).
We call $p'_i$ \emph{lower} if $p'_i < p'_{i+1}$ and \emph{upper} otherwise.
The cycle $C'$ is now homotopic to the cycle $C''$ 
formed by replacing each lower vertex $p'_i$ with an atom $a_i \leqs p'_i$.  (The homotopy is 
obtained by considering the simplices $a_i\leqs p_i < p_{i-1}$ and $a_i\leqs p_i<p_{i+1}$.)

Let $C'' = (a_1, q_1, a_2,\ldots a_n, q_n, a_1)$ (with the $a_i$ atoms).  Then
$C''$ is homotopic to the edge cycle formed by replacing each upper vertex $q_i$ by 
$\gen{a_i, a_{i+1}}$.  (The homotopy is obtained by considering the simplices
$a_i\leqs \gen{a_i, a_{i+1}}\leqs q_i$ and $a_{i+1}\leqs \gen{a_i, a_{i+1}}\leqs q_i$.)
$\hfill \Box$ 

\vspace{.15in}
We call edge cycles of the form described in the lemma \emph{atomic} cycles.

\vspace{.15in}
\noindent {\bf Alternate Proof of Corollary~\ref{2-generator}.}
Let $P$ be an atomized poset in which no three atoms generate.  Theorem~\ref{k-connected-atom}
shows that $P$ is connected, so it suffices to show that $\pi_1 (P) = 1$.
By the lemma any loop $L$ in $\abs{P}$ is homotopic to an atomic cycle 
$C = (a_1, \gen{a_1, a_2}, a_2,\gen{a_2, a_3},\ldots a_n, \gen{a_n, a_1}, a_1)$
in which the $a_i$ are atoms.  Note that $C$ has even length, and that if the length
of $C$ is two or four, $C$ is simply a path followed by its inverse and is thus nullhomotopic.

If the length of $C$ is greater than 4, then $C$ contains the path  
$\gamma = a_1<\gen{a_1, a_2}>a_2<\gen{a_2, a_3}>a_3$, and since $\gamma$
lies in the cone under $\gen{a_1, a_2, a_3}\neq P$ it is homotopic (relative to its endpoints) 
to any other path between
$a_1$ and $a_3$, provided the second path also lies under $\gen{a_1, a_2, a_3}$.
In particular, $\gamma$ is homotopic to the path $(a_1\leqs \gen{a_1, a_3}\geqs a_3)$.
Thus $C$ is homotopic to a shorter cycle, which
is still atomic, and repeating the process will eventually provide a null-homotopy of $C$.
$\hfill \Box$ 

\vspace{.15in}
We now present a result guaranteeing simple connectivity of an atomized poset
under weaker conditions than those of
Corollary~\ref{2-generator}.  

The following result is standard; see~\cite{Spanier-alg-top}.

\begin{lemma}$\label{presentation}$  Let $\Delta$ be a simplicial complex, and 
let $T\subset \Delta^{\leqs\, 1}$ be a maximal tree (i.e. a spanning tree).  

Then $\pi_1(\Delta)$ has a presentation with a generator for each (ordered) edge
$(u,v)$ with $\{u,v\}\in \Delta^{\leqs\, 1}$, and with the following relations:
\begin{enumerate}
	\item $(u,v) = 1$ if $\{u,v\}\in T$,
	\item $(u,v)(v,u) = 1$ if $\{u,v\}\in \Delta^{\leqs\, 1}$,
	\item $(u,v)(v,w)(w,u) = 1$ if $\{u,v,w\}\in \Delta^{\leqs\, 2}$.
\end{enumerate}
\end{lemma}

When the simplicial complex in question is $\M{G}$ for some non-cyclic group $G$,
we will always set the tree $T$ to be the collection of all edges $\{1,g\}$
($g\in G$).  (Note that since $G$ is non-cyclic, all such edges exist.)
We will refer to the resulting presentation for $\pi_1 (\M{G})\isom \pi_1 (\C{G})$
as the \emph{standard presentation}.

\begin{proposition}$\label{presentation-atom}$
Let $P$ be an atomized poset and assume that no two atoms generate $P$.
Say there exists an atom $a_0\in \A (P)$ with the following property:
\vspace{.05in}
\newline (1) \indent For any elements $a_1, a_2\in \A (P)$, there exists an 
element $a_3\in \A (P)$ 
such that $\gen{a_1, a_2, a_3}$, $\gen{a_0, a_1, a_3}$, $\gen{a_0, a_2, a_3}\neq P$.

\vspace{.05in}
Then $\abs{P}$ is simply connected.
\end{proposition}
{\bf Proof.}  First, note that in the 1-skeleton of $\M{P}$
all possible edges between vertices exist (and in particular $\M{P}$ is connected).
We wish to apply Lemma~\ref{presentation} to $\M{P}$, so we must choose a maximal tree
in $\M{P}^{\leqs\, 1}$.  Let $T$ be the star at
the vertex $a_0$, i.e. $T$ consists of all edges of the form $\{a_0,a\}$ with $a\in \A (P)$.  
The generators given in the lemma consist of ordered edges $(a_1, a_2)$ ($a_1, a_2\in \A (P)$)
and we will now show that each generator is trivial.

There are two cases.  First, if $\gen{a_0, a_1, a_2}\neq P$, then these vertices form a simplex 
in $\M{P}$.  So we obtain the relations $(a_1,a_2) = (a_1,a_0)(a_0,a_2)=1$. 

Next, say $\gen{a_0, a_1, a_2} = P$.  
Let $a_3$ be the element guaranteed by the
hypothesis.  Then, as above, the generators corresponding to the edges $\{a_1,a_3\}$ and
$\{a_2,a_3\}$ are trivial.  
Furthermore, since $\gen{a_1, a_2, a_3} \neq P$ these vertices form a simplex in $\M{P}$, and
we have the relations $(a_1,a_2) = (a_1,a_3)(a_3,a_2) = 1$. 

Thus all generators for $\pi_1 (\M{P})$ are trivial, and hence $\M{P}$ and $\abs{P}$ are simply
connected.
$\hfill \Box$

\vspace{.15in}
In the case of the coset poset, it is easy to check that Proposition~\ref{presentation-atom}
specializes to the following result (using the standard presentation, i.e. $a_0 = 1$).

\begin{corollary}$\label{2-gen-presentation}$
Let $G$ be a non-cyclic group in which the following condition holds: 
\vspace{0.05in}
\newline (2) \indent For any $x,y\in G$ such that 
$\gen{x,y} = G$, there exists $z\in G$ such that $\gen{z,x}$, $\gen{z,y}$, 
$\gen{z^{-1}x,z^{-1}y} \neq G$.  

\vspace{.05in}
Then $\C{G}$ is simply connected.
\end{corollary}

We now consider the cases in which these last two results apply.  Of course, if $P$ is not
generated by three atoms (respectively $G$ is not generated by two elements) then (1)
(respectively (2)) is satisfied, but Theorem~\ref{k-connected-atom} also applies.  
Additionally, if there is an element $a_0\in P$ such that $\gen{a_0, x, y}\neq P$ for any
$x,y\in P$, then ($\ast$) is satisfied.  (In the case of the coset poset, this reduces to saying
$G$ is not a 2-generator group.)

The simple group $A_5$ satisfies (2) and hence $\C{A_5}$ is simply connected.  The
details, which are not difficult, appear in~\cite{Ramras-thesis}.  In Section~\ref{simple-groups}
we will give another proof of simple connectivity, which eliminates some of the computation.

We note one final situation in which (2) is satisfied.  Say $G = H\cross K$, where
$K$ and $H$ are non-cyclic groups.  Then if $x = (h_1, k_1)$ and $y = (h_2, k_2)$, setting
$z = (h_1, k_2)$ we have $\gen{x,z}$, $\gen{y,z}$ and $\gen{z^{-1}x, z^{-1}y}\neq G$.  
In Section~\ref{extensions} we will go a bit further and classify direct products with simply
connected coset posets (Theorem~\ref{direct-products}).
$\hfill \Box$

\section{The Coset Poset of an Extension}$\label{extensions}$

In this section we will consider the coset poset of a
(non-simple) group $G$ in terms of an extension $N\normal G \to G/N$.

\begin{definition}
For any semi-direct product $G = H\rtimes K$ ($H$ and $K$ arbitrary groups),
let $f:G \to H$ be the function $f(h,k) = h$, and let $\pi: G\to K$
be the quotient map.  We call a coset $gT\in \C{G}$ \emph{saturating} if
$\pi (gT) = K$ and the only $K$-invariant subgroup of $H$ that contains $f(T)$ is $H$
itself (a subgroup $I\leqs H$ is called $K$-invariant if each element of $K$ induces
automorphisms of $I$).
\end{definition}

The direct product $P\cross Q$ of posets $P$ and $Q$
is the Cartesian product of $P$ and $Q$ together with the ordering 
$(p,q)\leqs(p',q')\iff p\leqs p'$ and $q\leqs q'$. 
The join $P*Q$ of $P$ and $Q$ is
the disjoint union of $P$ and $Q$ together with the ordering that induces the original 
orderings on $P$ and $Q$ and satisfies $p<q$ for each $p\in P$, $q\in Q$.  There
are canonical homeomorphisms $\abs{P\cross Q}\homeo \abs{P}\cross \abs{Q}$
and $\abs{P*Q}\homeo \abs{P}*\abs{Q}$, so long as one takes the associated compactly generated
topology on the right hand side when $P$ and $Q$ are not locally countable 
(see~\cite{Walker-canonical-homeomorphisms}).

\begin{lemma}$\label{semi-direct-products}$
Let $G = H\rtimes K$, with $H$ and $K$ non-trivial groups.
Let $\mC_0 (G)$ be the poset of all non-saturating cosets and let $\mC_K (H)$ 
denote the poset of all cosets of proper, $K$-invariant subgroups of $H$.  
Then $\mC_0 (G)$ is homotopy equivalent to $\mC_K (H)*\C{K}$.
\end{lemma}
{\bf Proof.}    Let
$\mC^+ (K)$ denote the set of all cosets in $K$ (including $K$ itself) and let
$\mC_K^+ (H) = \mC_K (H) \cup \{H\}$.  Then
if $\mC_{00} (H\rtimes K)$ denotes the set of all proper cosets of the form 
$(g,h) I\rtimes J$ (with $I\leqs H$ $K$-invariant and $J\leqs K$), one checks that the
map
$$(x,y) I\rtimes J \mapsto (xI,yJ)$$
is a well-defined poset isomorphism 
$$\mC_{00} (H\rtimes K) \stackrel{\isom}{\maps} \mC_K^+ (H) \cross \mC^+ (K) - \{(H, K)\}.$$
The latter is homeomorphic to $\abs{\mC_K (H) * \C{K}}$ 
(see~\cite[Proposition 1.9]{Quillen-p-subgroups} or the proof 
of~\cite[Prosition 2.5]{Kratzer-subgroups}).
Finally, we have an increasing poset map $\Phi$ from $\mC_0 (H\rtimes K)$ onto 
$\mC_{00} (H\rtimes K)$ given by
$$\Phi((x,y) T) = (x,y) \hat{f}(T)\rtimes \pi (T)$$
where $\hat{f} (T)$ is the smallest $K$-invariant subgroup containing $f(T)$ 
(i.e. the intersection of all invariant subgroups containing $f(T)$).
This map is a homotopy equivalence by~\cite[Corollary 10.12]{Bjorner-top-methods}.
$\hfill \Box$

\vspace{.15in}
In the case where the action of $K$ on $H$ is trivial, $f$ becomes the quotient
map $H\cross K\twoheadrightarrow H$ and all subgroups of $H$ are $K$-invariant.  
In this case $\mC_0 (H\cross K)$
is the poset of all cosets which do not surject onto both factors, and we obtain
the following result which appears as~\cite[Lemma 5, Proposition 12]{Brown-coset-poset}.

\begin{lemma}[Brown]$\label{brown-direct-products}$
For any finite groups $H$ and $K$, $\mC_0 (H\cross K)\heq \C{H}*\C{K}$.
If $H$ and $K$ are coprime (i.e. have no isomorphic quotients other than the trivial group)
then in fact $\C{H\cross K}\heq \C{H} * \C{K}$.
\end{lemma}

Coprimality implies that there are no saturating subgroups (and hence no saturating cosets);
see~\cite[Section 2.4]{Brown-coset-poset}.  

We will now use Lemma~\ref{semi-direct-products} to show that most semi-direct products
have simply connected coset posets.  First we need the following simple lemma, which 
appears (for finite groups) as~\cite[Proposition 14]{Brown-coset-poset}.  
The result extends, with the same proof, to infinite groups 
(the argument for $\Z$ requires a simple modification).

\begin{lemma}$\label{coset-poset-connectivity}$
Let $G$ be a group. Then $\C{G}$ is connected unless $G$ is cyclic of prime-power
order.
\end{lemma}

Given a semi-direct product $G = H\rtimes K$, we say that $H$ is $K$-cyclic
if there is an element $h\in H$ which is not contained in any proper $K$-invariant
subgroup of $H$.

\begin{proposition}$\label{pi-semi-direct}$
Let $G = H\rtimes K$ with $H$ and $K$ non-trivial groups.  If $K$ is not finite cyclic, or if 
$H$ is not $K$-cyclic, then $\C{G}$ is simply connected.
Furthermore, if $G$ is a torsion group and $K$ is not cyclic of prime-power order, $\C{G}$
is simply connected.
\end{proposition}
{\bf Proof.}  By Lemma~\ref{semi-direct-products}, we have 
$\mC_0 (G)\heq \mC_K (H) * \C{K}$.  We claim that in each of the above cases, $\mC_0 (G)$ is simply
connected.  The join of a connected space and a non-empty space is always simply connected
(see~\cite{Milnor-univ-bundles-2}), and both $\mC_K (H)$ and $\C{K}$ are always
non-empty ($\{1\}\in \mC_K (H)$) so it suffices to show that one or the other is connected.
If $H$ is not $K$-cyclic, then for every element $h\in H$ there is a $K$-invariant subgroup
$I_h < H$ containing $h$, and thus we have a path $hT\geqs \{h\} \leqs I_h\geqs \{1\}$ joining
any coset $hT\in \mC_K (H)$ to the trivial subgroup.  In the other cases, 
Lemma~\ref{coset-poset-connectivity} shows that $\C{K}$ is connected.  Thus $\pi_1 (\mC_0 (G)) = 1$
in each case.

We now show that every loop in $\C{G}$ is nullhomotopic.  By Lemma~\ref{cyclic-cycle} 
it suffices to consider edge cycles of the form 
$$C = (\{x_1\}, x_1 T_1, \{x_2\}, x_2 T_2, \ldots, \{x_n\}, x_n T_n, \{x_1\}),$$
where the $T_i$ are cyclic.
If $H$ is not $K$-cyclic, then every cyclic subgroup of $G$ lies in $\mC_0 (G)$ and
hence $C$ lies in $\abs{\mC_0 (G)}$ and must be nullhomotopic.  

Next, say $K$ is not a finite cyclic 
group.  If none of the vertices of $C$ are saturating cosets, we are done, so assume some coset 
$x_i T_i$ saturates.  Since
$T_i$ is cyclic and $K$ is not finite cyclic, we must have $T_i\isom K\isom \Z$.  
So no subgroup of $T_i$ surjects onto $K$, and hence no subcoset of $x_i T_i$ saturates. 
Thus adding $x_iT_i$ to $\mC_0 (G)$ cones off a copy of $\C{T_i}$, which 
(by Lemma~\ref{coset-poset-connectivity}) is connected.  Since
the union of two simply connected spaces with connected intersection is simply connected,
we see that $\mC_0 (G)\cup \{x_i T_i\}$ is simply connected 
($\ord{
       \mC_0 (G)
          \cup 
       \{x_i T_i\}
      } 
 = \ord{
        \mC_0 (G)
       } 
      \cup
   \ord{ 
        \C{G}_{\leqs x_i T_i}
       }
$).
Repeating the process eventually shows that $C$ lies in a simply connected poset, and hence
is nullhomotopic.

If $K$ is not cyclic of prime-power order and $G$ is a torsion group, each $T_i$ is 
finite.  Thus there are finitely many cosets in the set
$S = \{xT: xT\subset x_i T_i {\rm\,\, for\,\, some\,\,} i\} - \mC_0 (G)$.  
If $xT$ is a minimum element of $S$, then adding $xT$ to $\mC_0 (G)$ cones of a copy of
$\C{T}$, which is connected ($T$ surjects onto $K$, which is not cyclic of prime-power order).  
This process may be repeated until we have added all
of $S$, and hence the poset $S\cup \mC_0 (G)$ is simply connected.  The cycle $C$ lies in this
poset and is thus nullhomotopic.
$\hfill \Box$

\begin{theorem}$\label{direct-products}$  Let $H$ and $K$ be non-trivial groups.
Then $\pi_1 (\C{H\cross K})\neq 1$ if and only if 
both $H$ and $K$ are cyclic of prime-power order.  
\end{theorem}
{\bf Proof.}  If both groups are cyclic of prime-power order, it is easy to 
check that there are just two complemented factors in any chief series for $H\cross K$
(in the sense of Theorem~\ref{solvable-coset-poset}) and the desired result then 
follows from that theorem.
In the other direction, apply Proposition~\ref{pi-semi-direct}.
$\hfill \Box$

\vspace{.15in}
The question of simple connectivity for the coset poset of a 
finite (non-trivial) semi-direct product is now
reduced to the case of products $H\rtimes \Z/p^n$, where $p$ is prime and $H$ is
$(\Z/p^n)$-cyclic.  When $H$ is solvable, Theorem~\ref{solvable-coset-poset} applies,
so we are most interested in the case where $H$ is a non-solvable group.  The simplest
example, then, is $S_5\isom A_5\rtimes \Z/2$.  In this case the coset poset is still
simply connected.  For the proof we will need the following result of 
Brown~\cite[Proposition 10]{Brown-coset-poset}, which we note extends (with the same proof)
to infinite groups.

\begin{lemma}[Brown]$\label{general-extension}$
For any group $G$ with normal subgroup $N$, there is a homotopy equivalence
$\C{G}\heq \C{G/N}*\C{G,N}$, where the latter poset is the collection of all
cosets $xH\in \C{G}$ which surject onto $G/N$ under the quotient map.
\end{lemma}

\begin{proposition}$\label{S_n}$ For $n>3$, the coset poset of $S_n$ is simply connected.
\end{proposition}
{\bf Proof.}  We have $A_n\normal S_n$ with $S_n/A_n\isom \Z/2$, 
so Lemma~\ref{general-extension} gives 
$\C{S_n}\heq \C{\Z/2}*\C{S_n,A_n}\homeo$ Susp $(\C{S_n, A_n})$, so it will
suffice to show that $\C{S_n,A_n}$ is connected. 

Since $A_n$ has index two in $S_n$, the elements of $\C{S_n,A_n}$ 
are exactly the cosets $xH$ where $H$ is not contained in $A_n$.  
Letting $S$ denote the set of proper subgroups of $S_n$
not contained in $A_n$, we have $\C{S_n, A_n} = \bigcup_{x\in S_n} \{xH:H\in S\}$.  

First we show that for each $x$ the poset $\{xH:H\in S\}$ is connected.  Each such poset
is isomorphic to $S$, so it suffices to consider this case.
Consider the action of $S_n$ on the set $\{1,\ldots, n\}$.  For each $i$, Stab$(i)\isom S_{n-1}$,
and for any $i,j$, Stab$(i)\cap$Stab$(j)\isom S_{n-2}$.  
All of these groups are in $S$ (since $n>3$), and hence for each $i, j$
there is a path in $S$ between Stab$(i)$ and Stab$(j)$.

We now show that for each $K\in S$ there exists $i$ such that we have a path 
(in $S$) from $K$ to Stab$(i)$.  
Since $K$ is not contained in $A_n$, there is an element $k\in K$ such that $k\notin A_n$.  First,
say $k$ is not an $n$-cycle.  Then if the
orbits of $\gen{k}$ have orders $a_1,\ldots, a_m$, there is a path of the form
$$K\geqs \gen{k}\leqs S_{a_1}\cross \cdots \cross S_{a_m} \geqs S_{a_1} \leqs S_{n-1}$$
connecting $K$ to the stabilizer of some point.  
If $k$ is an $n$-cycle, we consider
two cases, depending on the parity of $n/2$ (note that $n$ must be even since the $n$-cycle
$k$ is not in $A_n$).  If $n/2$ is odd, then $k^2$ is a product of two odd-length cycles
and hence is not in $A_n$.  Replacing $k$ by $k^2$ reduces to the above case.  When
$n/2$ is even, assume without loss of generality that $k = (1\,\,\,\, 2\cdots n)$.
The reader may check that (setting $m=n/2$) the element 
$t = (m\,\,\,\, m+2)(m-1\,\,\,\, m+3)\cdots (2\,\,\,\, 2m)$
lies in the normalizer of $\gen{(1\,\,\,\, 2 \cdots n)}$.  This element is a product of
$m-1$ disjoint transpositions, and hence is not in $A_n$.  Thus we have a path of the form
$$K\geqs \gen{k}\leqs N (\gen{k})\geqs \gen{t}\leqs S_2\cross S_2\cross\cdots\cross S_2
\geqs S_2\leqs S_{n-1}.$$

To complete the proof, we must find a path in $\C{S_n, A_n}$ joining $C_x = \{xH:H\in S\}$ to
$C_y = \{yH:H\in S\}$ for each $x,y\in S_n$.  If $x^{-1}y\notin A_n$ then $x\gen{x^{-1}y}$
lies in $C_x \cap C_y$ and we are done.  If $x^{-1}y\in A_n$, then either $x, y\in A_n$
or $x, y\notin A_n$.  In the first case, take $z\notin A_n$.  Then $x^{-1}z, z^{-1}y\notin A_n$
and hence $C_x\cap C_z \neq \emptyset$ and $C_z\cap C_y\neq \emptyset$.  Since $C_z$ is connected,
this yields a path joining $C_x$ to $C_y$.  If $x,y\notin A_n$, then choosing $z\in A_n$ we
may complete the proof in a similar manner.
$\hfill \Box$

\begin{remark} Analogous arguments shows that $\C{\Z}$ and $\C{\Z\rtimes \Z/2}$, where $\Z/2$ acts 
by inversion, are simply connected (again, the normal subgroup of index two is key).  
\end{remark}

Turning to the case of a general group extension, we note a simple consequence of
Lemma~\ref{general-extension}. 

\begin{proposition}$\label{coset-quotient}$
Let $G$ be a group with quotient $\overline{G}$, and assume that
$\C{\overline{G}}$ is $k$-connected.  Then $\C{G}$ is $k$-connected as well.
\end{proposition}
{\bf Proof.}  This follows immediately from the lemma, noting that the join of a 
$k$-connected space with any space is still $k$-connected~\cite{Milnor-univ-bundles-2}.
$\hfill \Box$

\section{Mayer-Vietoris Sequences and $H_* (\C{G})$}
$\label{m-v}$

We will now present some results based on Mayer-Vietoris sequences for the homology of 
$\C{G}$.  In particular, we will obtain lower bounds on the rank of $H_1 (\C{G})$
for certain groups $G$.  The results in this section are most interesting
for infinite groups, as the classes of finite groups to which the results apply are rather small.

From here on, we let $H_n(X)$ denote the simplicial homology of the space $X$ with coefficients 
in $\Z$, and we let $\tilde{H}_n (X)$ denote reduced simplicial homology 
(again with coefficients in $\Z$).

\begin{theorem}$\label{M-V1}$
Let $G$ be a non-cyclic group with a cyclic maximal subgroup $M$
of prime-power order, and
let $o(M) = p^n$.  Then $H_1(\C{G})$ has rank at least $(p-1) (G:M)$ (or infinite rank
if $G$ is infinite), and in particular 
$\C{G}$ is not simply connected.
\end{theorem}
{\bf Proof.}  Let $\{x_i\}_{i\in I}$ be a set of left coset representatives for $M$.
Let $X = \ord{\C{G}}$, $Y = \ord{\C{G} - \{x_i M\}_{i\in I}}$, 
and $Z = \ord{\cup_{i\in I} \C{G}_{\leqs x_i M}}$.
Then we have $X = Y\cup Z$, yielding a Mayer-Vietoris sequence
$$\cdots \maps \tilde{H}_1(X) 
         \stackrel{\partial}{\maps} \tilde{H}_0 (Y\cap Z) 
         \stackrel{h}{\maps} \tilde{H}_0 (Y) \oplus \tilde{H}_0 (Z) 
         \maps \tilde{H}_0 (X)$$
on (reduced) homology groups.  

By Lemma~\ref{coset-poset-connectivity}, $X$ is connected, so $\tilde{H}_0 (X) = 0$ and
$h$ is a surjection.
Also, $Z$ is a union of $(G:M)$ disjoint cones, so 
$\tilde{H}_0 (Z) \homeo \Z^{(G:M) - 1}$ and
we claim that $Y$ is connected (so that $\tilde{H}_0(Y) = 0$).  

Now, 
$$Y\cap Z = \ord{\cup_{i\in I} \C{G}_{<x_i M}} \homeo \coprod_{i\in I} \ord{\C{M}},$$
where $\coprod$ denotes the disjoint union.  
First, assume that $o(G)<\infty$.  Then the coset
poset of $M\homeo \Z/p^n$ has $p$ contractible components 
(the cones under the cosets of the unique maximal subgroup) and
thus $\tilde{H}_0 (Y\cap Z) \homeo \Z^{p (G:M) - 1}$.
Substituting these
values into the above sequence yields an exact sequence
$$\cdots \maps \tilde{H}_1(X) \stackrel{\partial}{\maps} \Z^{p(G:M) - 1} 
         \stackrel{h}{\maps} \Z^{(G:M) - 1} \maps 0.$$
So Im$\, \partial = $ Ker$\, h \homeo \Z^{p(G:M) - 1 - ((G:M) - 1)} = \Z^{(p-1)(G:M)}$ 
and we see that the rank of $H_1(X)$ is at least $(p-1)(G:M)$, as desired.  

If $G$ is infinite, Ker$\, h$ contains all homology classes of the form
$[xM' - y M']$ where $M'<M$ is the unique maximal subgroup of $M$ and $x\equiv y$ (mod $M$)
(the image of $[xM' - yM']$ in $\tilde{H}_0(Y)$ is trivial because, as we will show, 
$Y$ is connected,
and its image in $Z$ is trivial because $\partial((yM', xM) + (xM, xM'))= xM' - yM')$, where
$\partial$ denotes the simplicial boundary).  Thus
Im $\partial = $Ker$\, h$ contains a free abelian group of infinite rank, and hence
the rank of $H_1 (X)$ is infinite as well (it is worth noting that maximality of the cyclic
group $M$ implies that $G$ is a two generator group and hence is countable).

To complete the proof we must show that $Y$ is connected.
Consider a coset $xH$, where
$H\neq M$.  If $\gen{x}\neq M$, then we have a path 
$$xH\geqs \{x\} \leqs \gen{x} \geqs \{1\}$$
connecting $xH$ to the identity.  Next, if $\gen{x} = M$ then choose some $g\in G$, $g\notin M$.
We now have a path
$$xH\geqs \{x\} \leqs x \gen{g} \geqs \{xg\} \leqs \gen{xg} \geqs \{1\}$$
in $Y$ connecting $xH$ to the identity.
$\hfill \Box$

\vspace{.15in}
The above result does not extend to cyclic groups,   
since the subcomplex $Y$ is no longer connected.  Nevertheless, Theorem~\ref{solvable-coset-poset}
can be used to show that 
the coset poset of a cyclic group $G$ with a maximal subgroup isomorphic to $\Z/p^n$ is still
not simply connected.  

There are many infinite groups to which Theorem~\ref{M-V1} applies.  Such 
groups are described in~\cite[Chapter 9]{Olshanskii-defining-relations}.  In
particular, for each sufficiently large prime $p$ there is a continuum of pairwise
non-isomorphic groups in which each non-trivial proper subgroup has order $p$. 

As mentioned at the beginning of the section, there are very few finite groups to which
the theorem applies.  The $p$-groups with a cyclic maximal subgroup have been
classified (see~\cite{Brown-cohomology}), and there are only a few types.  In 
fact, any other finite group $G$ with a cyclic maximal subgroup $M$ of prime-power 
order is either a semi-direct product $A\rtimes \Z/p^n$, where $A$ is elementary abelian, or of
the form $\Z/p^n \rtimes \Z/q$ where $q$ is prime.  This can be shown as follows.  
If $M\normal G$, then clearly $G\isom \Z/p^n \rtimes \Z/q$.  Otherwise, 
Herstein~\cite{Herstein-remark} has shown (by very elegant and elementary
methods) that a group with an abelian maximal subgroup is solvable.  Now, $M$ is a Sylow 
$p$-subgroup of $G$ and we have $N_G (M) = M$, so $M$ lies in the center of its normalizer
and must have a complement (by Burnside's Theorem~\cite[p. 289]{Robinson-theory-of-groups}).  
So $G= T\rtimes M$ for some $T<G$, and maximality of $M$ implies that no subgroup of $T$ is 
invariant under the action of $M$, i.e. $T$ is a minimal normal subgroup.  Since $G$ is solvable, 
this implies that $T$ is in fact elementary abelian.  The simplest interesting example of such
a group is $A_4\isom (\Z/2)^2\rtimes \Z/3$, and other examples may be constructed by taking an
appropriate generator $x$ of the field extension $\F_q\subset \F_{q^n}$ and letting $x$ act by multiplication
on the additive group of $\F_{q^n}$.

The interested reader may now use Theorem~\ref{solvable-coset-poset} to compute 
the exact homotopy type of $\C{G}$ for any finite group $G$ to which the theorem applies 
(it is worth noting that the bound on the rank of $H_1$ is not very good).  The theorem is still
useful in this endeavor, though, since it proves that there is exactly one complemented factor
in the chief series for $G$ and hence simplifies the computation of the homotopy type of $\C{G}$.

\begin{problem} Find a finite non-cyclic group $G$ with $\C{G}$ not simply connected
and such that no maximal subgroup of $G$ is cyclic of prime-power order.  A non-solvable
example would be of particular interest.
\end{problem}

The cyclic groups of order $p^n q^m$ ($p,q$ prime) have non-simply connected
coset posets (Theorem~\ref{direct-products}) but have no maximal subgroup which is cyclic
of prime-power order (unless $n=m=1$).

The following lemma shows a (somewhat weak) connection between the coset poset and
the subgroup poset of a finite group.  Further relationships are discussed in the final section.

\begin{lemma}$\label{homology-surjection}$
Let $G$ be any group, and for $g\in G$, let $\C{G}_g$ denote
the poset $\C{G} - \{g\}$.

If there exists $g\in G$ with $\tilde{H}_n (\C{G}_g) = 0$,
then there is a surjection 
$$\tilde{H}_{n+1} (\C{G}) \twoheadrightarrow \tilde{H}_n (L(G)).$$
\end{lemma}

Note that for any $g, h\in G$, the posets $\C{G}_g$ and $\C{G}_h$ are isomorphic
(consider the map given by left-multiplication by $hg^{-1}$).  In particular,
$\C{G}_g \isom \C{G}_1$.

\vspace{.15in}
{\bf Proof.}  Let $X = \ord{\C{G}}$, $Y = \ord{\C{G}_{\geqs \{g\}}}$, and 
$Z =  \ord{\C{G}_g}$.  
Note that $Y$ is contractible (so $\tilde{H}_n (Y) = 0$) and 
that $Y\cap Z \isom L(G)$.

Since $X = Y\cup Z$, we get a Mayer-Vietoris sequence
$$\cdots \maps \tilde{H}_{n+1} (X) \stackrel{\partial}{\maps} \tilde{H}_n (Y\cap Z) 
     \maps \tilde{H}_n (Y)\oplus \tilde{H}_n (Z) \maps\cdots
$$
and since $\tilde{H}_n (Y) = \tilde{H}_n (Z) = 0$, $\partial$
is a surjection.
$\hfill \Box$

\begin{theorem}$\label{M-V2}$
Let $G$ be a group and let $n$ be the cardinality of the set of path-components of $L(G)$.
Then $H_1(\C{G})$ has rank at least $n-1$ (or infinite rank if $n$ is infinite).  
In particular, if $L(G)$ is disconnected then $\C{G}$ is not simply connected.
\end{theorem}
{\bf Proof.}  If $G$ is cyclic then $L(G)$ is connected (unless $G\isom \Z/pq$, $p,q$ prime) 
and the result is trivial.  When $G\isom \Z/pq$, the result follows
from Theorem~\ref{solvable-coset-poset}.

We now assume that $G$ is not cyclic.  By the lemma, it will suffice to prove that
$\C{G}_1$ is connected.    
Choose some element $x\in G$, $x\neq 1$.  Then any vertex $yH$
in $\C{G}_1$ can be connected to $\{x\}$ via the path 
$$yH \geqs \{y\} \leqs x\gen{x^{-1} y}\geqs \{x\}$$
($yH\neq \{1\}$ so we may assume $y\neq 1$).  
$\hfill \Box$

\vspace{.15in}
For finite groups, it appears that this result is eclipsed by Theorem~\ref{M-V1}.  See
Lemma~\ref{subgroup-connectedness} and discussion following it for details.

There is again an interesting class of infinite groups to which Theorem~\ref{M-V2} applies
(in addition to the groups described above, to which Theorem~\ref{M-V1} also applies).
In~\cite[Chapter 9]{Olshanskii-defining-relations}, it is shown there that there is a continuum
of non-isomorphic infinite groups $G$ 
all of whose non-trivial proper subgroups are infinite cyclic.
Furthermore, in each of these groups, any two maximal subgroups intersect trivially 
and hence $L(G)$ is disconnected in every case.
\section{2-transitive Covers and Simple Groups}$\label{simple-groups}$

In this section, we will examine the homotopy type of the coset posets of the finite simple 
groups $PSL_2 (\F_5) \isom A_5$ and $\psls$, utilizing the notion of a 2-transitive cover.
Some computational details will be omitted, and a complete presentation of these results (and
of all the necessary background information on $\psls$) may be found in~\cite{Ramras-thesis}.

\subsection{2-transitive Covers}
\begin{definition} Let $G$ be a group.  We call a collection of subgroups
$S\subset L(G)$ a \emph{cover} of $G$ if every element of $G$ lies in some subgroup $H\in S$.  
We call $S$ \emph{2-transitive} if for each $H\in S$, the action of $G$ on the left cosets of $H$
is 2-transitive.

In addition, we say that a 2-transitive cover $S$ is \emph{$n$-regular} 
if for each $H\in S$ there 
is an element $g\in G$ with $o(g)=n$ whose action on $G/H$ is non-trivial (this is equivalent 
to requiring that no subgroup $H\in S$ contains all the elements in $G$ of order $n$).
\end{definition}

Recall that the \emph{standard presentation} for
$\pi_1 (\M{G})\isom \pi_1 (\C{G})$ is the presentation obtained from Lemma~\ref{presentation}
using as maximal tree the star at the vertex $1$ (assuming $G$ is non-cyclic).

\begin{definition} Let $G$ be a non-cyclic group, and say $G$ contains elements of order $n$.  
We say that $\M{G}$ is \emph{$n$-locally simply connected} if each generator $(g,h)$
(in the standard presentation for $\pi_1 (\M{G})$) with $o(g) = n$ is trivial.
\end{definition}

\begin{proposition}$\label{covers}$
Let $G$ be a non-cyclic group containing elements of order $n$.  
If $\M{G}$ is $n$-locally simply connected and $G$ admits an $n$-regular 2-transitive cover,
then $\M{G}$ (and hence $\C{G}$) is simply connected.
\end{proposition}
{\bf Proof.}
Let $\{g_1,g_2\}$ be any edge of $\M{G}$.  
We must show that the corresponding generators in the standard presentation for
$\pi_1 (\M{A_5})$ are trivial.
It will suffice to show that there exists an element
$z$ of order $n$ and a subgroup $K\in L(G)$ such that $g_1\equiv g_2 \equiv z$ (mod $K$).  
(Then the set $\{g_1,g_2,z\}$ forms a 2-simplex in $\M{G}$ and we have
$(g_1,g_2) = (g_1,z)(z,g_2) = 1$ since $\M{G}$ is $n$-locally simply connected.)

Let $S$ be an $n$-regular 2-transitive cover of $G$.  Then there exists $H\in S$ with
$g_1^{-1}g_2\in H$, and there is an element $z\in G$ with $o(z) = n$ which acts non-trivially
on the left cosets of $H$.  Let $H, x_1 H, \ldots, x_k H$ denote these cosets.
Since $z$ acts non-trivially on $G/H$ and $G$ acts 2-transitively,
some conjugate of $z$ sends $H$ to $x_i H$ ($1\leqs i\leqs k$), i.e. there is an element of 
order $n$ in every set $\{g\in G: g\cdot H = x_i H\} = x_i H$ ($1\leqs i\leqs k$).  
So we have found an element $z^{\alpha}\in g_1H = g_2H$, as desired.
$\hfill \Box$

\subsection{The Coset Poset of $A_5$}$\label{C(A_5)}$

We will establish the following:

\begin{proposition}$\label{A_5}$
The coset poset of $A_5\isom PSL_2 (\F_5)$ has the homotopy type of a bouquet of 
1560 two-dimensional spheres.
\end{proposition}

It can 
be checked that there are 1018 proper cosets in $A_5$, and hence $\C{A_5}$ has
1018 vertices.  Evidently, $\C{A_5}$ is far too large to admit direct analysis.  

Shareshian [unpublished manuscript]
has provided a proof of this result using the theory of shellability.  
The proof given here is somewhat simpler than Shareshian's argument.

Following Shareshian, we will show that 
$\mC = \C{A_5}$ has the homotopy type of a two-dimensional complex.  For
this portion of the proof we will work directly with $\mC$.  We will
show that $\mC$ is simply connected by applying Proposition~\ref{covers}.  
To show that $\C{A_5}$ has the homotopy type of a bouquet of 2-spheres,
we appeal to the general fact that a $k$-dimensional complex which is $(k-1)$-connected
is homotopy equivalent to a bouquet of $k$-spheres.  The number of spheres in the bouquet 
can be calculated from the Euler characteristic 
$\chi (\C{A_5})$, computed in~\cite{Brown-coset-poset}.

As mentioned in Section~\ref{connectivity}, the reader may also check that $A_5$
satisfies the conditions of Corollary~\ref{2-gen-presentation}.  This provides
a second proof of simple connectivity.  The proof given below eliminates some
of the computation.

\begin{claim}$\label{A_5-dimension}$
Let $\mC^-$ denote the poset $\mC$ with all cosets of all copies of $D_4$
removed.  Then $\ord{\mC^-}$ is two-dimensional and $\mC^- \heq \mC$.
\end{claim}
{\bf Proof (Shareshian).}  Quillen's Theorem A (see~\cite[Theorem 1.6]{Quillen-p-subgroups} 
or~\cite[Theorem 10.5]{Bjorner-top-methods})
shows that the inclusion $\mC^- \injects \mC$ is a homotopy equivalence 
(there is a unique subgroup lying above any copy of $D_4$, namely a copy of $A_4$), 
so it remains to check that $\mC^-$ is two-dimensional.  This follows easily from the fact that
each maximal subgroup of $A_5$ is isomorphic $A_4, D_{10}\isom Z/5\rtimes \Z/2$, or $S_3$. 
$\hfill \Box$

\vspace{.15in}
\begin{claim}$\label{pi-A_5}$
The coset poset of $A_5$ is simply connected.
\end{claim}
{\bf Proof.}
We will show that $A_5$ admits a 2-regular 2-transitive cover.  The interested
reader may check that $\M{A_5}$ is 2-locally simply connected using the method in
the proof of Proposition~\ref{presentation-atom}.  Up to automorphism, there are only 
several cases to check.

We claim that the set $\{H\in L(A_5): H\isom A_4 {\rm \,\, or\,\,} H\isom D_{10}\}$
is a 2-regular 2-transitive cover of $A_5$.  Each copy of $A_4$
is the stabilizer of a point in $\{1,2,3,4,5\}$, and the action of $A_5$ on this set is
2-transitive.  Each copy of $D_{10}$ is the stabilizer of a point under the action of 
$PSL_2(\F_5)\isom A_5$ on $\F_5\cup \infty$, and this action is 2-transitive as well
(\cite{Burnside-theory-of-groups},~\cite{Ramras-thesis}).  
In each case, all non-trivial elements act non-trivially, and every element in $A_5$
lies either in a $D_{10}$ or in an $A_4$.
$\hfill \Box$

\subsection{The Coset Poset of $\psls$}$\label{psl7}$

We now consider the finite simple group $G = PSL_2 (\F_7)$, and 
show that $\C{G}$ is simply connected.  Other facts about the homotopy type of $\C{G}$ are 
discussed at the end of the section.

For basic facts about the groups $\psl$, we refer 
to~\cite{Burnside-theory-of-groups} 
(see also~\cite{Ramras-thesis} and~\cite{Suzuki-group-theory-I}).  
We write the elements of $\psl$ as 
``M\"{o}bius transformations" $f:\prj\to \prj$,
$x\mapsto \frac{ax+b}{cx+d}$, where $a,b,c,d\in \F_p$, det$(f) = ad-bc = 1$, and
$\infty$ is dealt with in the usual manner.  The action of $\psl$ on $\prj$
is 2-transitive.  For a proof of the following result, see~\cite{Burnside-theory-of-groups}
or~\cite{Ramras-thesis}.

\begin{lemma}$\label{maximal-subgroups}$
Any maximal subgroup of $G$ is either the stabilizer of a point in $\F_7 \cup \{\infty\}$
(and is isomorphic to $\Z/7\rtimes \Z/3$) or is isomorphic to $S_4$.  The two conjugacy
classes of subgroups isomorphic to $S_4$ are interchanged by the
transpose-inverse automorphism of $GL_3 (\F_2)\isom G$.
\end{lemma}

The following lemma will help to minimize the amount of computation in the proof of simple 
connectivity.

\begin{lemma}$\label{psls-orders}$
Let $\alpha = \frac{ax+b}{cx+d}$ be any non-trivial element of $\psls$.  
We define the trace-squared of $\alpha$ to be $tr^2 (\alpha) = (a+d)^2$ 
(note that this is the square of the
trace of either representative of $\alpha$ in $SL_2 (\F_7)$, and thus is well-defined).
The order of $\alpha$ is then determined as follows:
$$o(\alpha) = \left\{ \begin{array}{ll}
			2, & tr^2 (\alpha) = 0 \\
			3, & tr^2 (\alpha) = 1 \\
			4, & tr^2 (\alpha) = 2 \\
			7, & tr^2 (\alpha) = 4. \\
		\end{array}
	\right.
$$
\end{lemma}
{\bf Proof.}  If $\alpha = \frac{ax+b}{cx+d}\in G$, we define
$disc(\alpha)$ to be the discriminant of the quadratic polynomial 
$(cx+d)x-(ax+b) = cx^2 + (d-a)x -b$ determined
by the equation $\frac{ax+b}{cx+d} = x$,
so that $disc(\alpha) = (d-a)^2 - 4 (-b)(c) = tr^2(\alpha) - 4$.

We consider only the elements in $G - $Stab$(\infty)$.  It is easy to check
the result on the remaining elements.  The elements of order seven in $G$ are exactly those
which fix one point in $\F_7\cup \infty$, and the elements of order three in $G$
are exactly those which fix two points in $\F_7\cup \infty$ 
(\cite{Burnside-theory-of-groups},~\cite{Ramras-thesis}).
Thus the elements of order seven in $G - $Stab$(\infty)$ are exactly those with 
$disc (\alpha) = 0$ and $tr^2 (\alpha) = 4$ and the elements of order
three in $G - $Stab$(\infty)$ are those with $disc(\alpha) \in \left(\F_7^*\right)^2$.  
Thus $o(\alpha) = 3 \iff disc(\alpha)\in \{1,2,4\} 
\iff tr^2 (\alpha)\in \{5, 6, 1\} \iff tr^2 (\alpha) = 1$.

Next, let $A\in SL_2 (\F_7)$ be a matrix representing $\alpha\in G$, and
let $\lambda_1, \lambda_2\in \bar{\F}_7$ be the eigenvalues of $A$
(where $\bar{\F}_7$ denotes the algebraic closure).  Then $o(\alpha) = 2
\Longrightarrow \lambda_1^2 = \lambda_2^2 = \pm 1$.  If $\lambda_1 = \lambda_2$ then
$A = \lambda_1 I$ and $\alpha = 1$, so $o(\alpha) = 2 \iff \lambda_1 = - \lambda_2
\iff tr(A) = tr^2 (\alpha) = 0$.
A similar but more lengthy calculation shows that 
$o(\alpha) = 4 \iff tr^2 (\alpha) = 2$.
$\hfill \Box$

\begin{proposition}$\label{psls-simply-connected}$
The coset poset of $PSL_2 (\F_7)$ is simply connected.
\end{proposition}
{\bf Proof.}  As in the case of $A_5\isom PSL_2 (\F_5)$, we will apply 
Proposition~\ref{covers}.  First we indicate the proof that $\M{G}$ is
2-locally simply connected.  We must check that each generator $(g,h)$
($o(g) = 2$) in the standard presentation for $\pi_1 (\M{G})$ is trivial.

Recall that $(g,h) =1$ if $\gen{g,h}\neq G$.
Given a group $H$, we define an automorphism class of generating pairs to be a set of 
the form $\{(\phi (x), \phi (y)): \phi \in Aut(H) {\rm\,\, and\,\,} \gen{x,y} = H\}$.
To show that each generator $(g,h)$ with $o(g) = 2$ is trivial, it suffices to check
one representative from each automorphism class of generators.
Letting $\Phi_{a,b}$ denote the number of automorphism classes of generators of $G$
in which all representatives $(x,y)$ satisfies $o(x) = a$ and $o(y) = b$,
M\"{o}bius inversion allows one to 
calculate $\Phi_{a,b}$ using the M\"{o}bius function of $G$ 
(as calculated in~\cite{Hall-eulerian}).  For these computations it is also necessary 
to know that $Aut(G)\isom PGL_2(\F_7)$ has order 336.  The method of calculation is described
in $\S1$ and $\S3$ of~\cite{Hall-eulerian}.

\vspace{.1in}
{\bf o(h) = 3:}  M\"{o}bius inversion shows $\Phi_{2,3} = 1$, 
and $g = \frac{x-2}{x-1}$, $h = \frac{4x}{2}$
represents the unique automorphism class of generators because $o(gh)=7$ 
and no subgroup contains elements of orders two and seven.  
Letting $z = \frac{3x-2}{-2x-3}$ we have $o(z) = 2$ and $g(-1) = h(-1) = z(-1) = -2$, so 
$\gen{z^{-1}g,z^{-1}h}\leqs Stab(-1)$.  
Since two elements of order two cannot generate a simple group,
$\gen{g,z}\neq G$, and since $o(hz) = 3$, the pair $(h,z)$ does not fall into the unique
automorphism class of generators with orders two and three.  Thus $\{z,g,h\}$ is a two-simplex
in $\M{G}$ and we have the relations $(g,h) = (g,z)(z,h) = 1$.

\vspace{.1in}
{\bf o(h) = 4:} M\"{o}bius inversion shows that $\Phi_{2,4} (G) = 1$.
Let $g = \frac{-1}{x}$, $h = \frac{4x+1}{-x}$ and $z = \frac{-2}{4x}$.  The argument is
now similar to the previous case (note that $g(0) = h(0) = z(0) = \infty$).

\vspace{.1in}
{\bf o(h) = 7:} This time we find that there are three automorphism
classes of generating pairs.  (Note that any pair of elements with these orders generates $G$.)
Let $h = x+1$ ($o(h) = 7$), and let $g = \frac{b}{cx}$ ($o(g) = 2$).
We have $hg = \frac{b}{cx} + 1 = \frac{cx+b}{cx}$
so $tr^2 (hg) = c^2$, which implies that the pairs
$(g_1 = \frac{-1}{x},\, h)$, $(g_2 = \frac{2}{3x},\, h)$, and 
$(g_3 = \frac{3}{2x},\, h)$
represent the three generating automorphism classes. 

Next, say there exist an elements $z_i\in$ Stab$(\infty)$ such that $o(z_i) = 3$ and
$\{h, g_i, z_i\}$ forms a simplex in $\M{G}$ ($i = 1,2,3$).  Then we have 
$\gen{h,z_i}\leqs$ Stab$(\infty)$ and $(g_i,z_i) = 1$ because $o(g_i) = 2$ and $o(z_i) = 3$.
So $(g_i,h) = (g_i,z_i)(z_i,h) = 1$.  We will now find such elements $z_i$.
Consider the equations $h(x) = g_i (x)$, i.e. $x+1 = \frac{b}{-b^{-1} x}$ ($b = -1, 2, 3$).  
Equivalently, (note that $\infty$ can never be a solution) we want to solve $x^2 + x + b^2 = 0$, 
and by examining the discriminant we see that solutions $x_i\in \F_7$ exist when $b = -1$ or $3$, 
but not when $b = 2$.
For $i=1,3$ there is an element $z_i\in$ Stab$(\infty)$ with $o(z_i) = 3$
and $z_i(x_i) = g_i(x_i) = h(x_i)$ ($z_i = \frac{2x + (2x_i + 4)}{4}$).
Thus $(g_1,h) = (g_3,h) = 1$.  

Finally, we must show that $(g_2, h) = (\frac{2}{3x}, x+1) = 1$.  
Letting $z_2 = \frac{2x}{4}$, we
have $o(z_2) = 3$ and $z_2\in$ Stab$(\infty)$, so it suffices to show that these three
elements lie in a proper coset, i.e. that $\gen{z_2^{-1} g_2, z_2^{-1} h}\neq G$.
We have $o(z_2^{-1} g_2) = 2$, $o(z_2^{-1} h) = 3$ and $o(z_2^{-1} g_2 z_2^{-1} h) = 4$,
so these elements do not lie in the unique automorphism class found at the start of the proof.

\vspace{.1in}
To complete the proof we will check that the set of maximal subgroups of $G$ is a 
2-regular 2-transitive cover.  By Lemma~\ref{maximal-subgroups}, the maximal subgroups are 
the stabilizers under the action of $G$ on $\F_7\cup \{\infty\}$, and the copies of $S_4$.
The action of $G$ on $\F_7\cup \{\infty\}$ is 2-transitive 
and all non-trivial elements act non-trivially.  
Next, consider the group $GL_3 (\F_2)\isom G$.  
It is easy to check that the stabilizer of a vector $v\in \F_2^3$ is isomorphic to $S_4$, 
and the action of $GL_3(\F_2)$ on $\F_2^3$ is 2-transitive.
Additionally, all non-trivial elements act non-trivially on $\F_2^3$.  
Since all copies of $S_4$ are conjugate under $Aut(G)$ (Lemma~\ref{maximal-subgroups}), 
this completes the proof.
$\hfill \Box$

\vspace{.15in}
There are two main difficulties in extending these results to 
$\psl$ for primes $p>7$.  It is not clear how to generalize the
ad hoc portion of the proof, in which we showed that $\M{G}$ is 2-locally simply connected.  
Also, for large $p$ the only 2-transitive action of $\psl$ is its standard action on $\prj$.  
This follows from~\cite[Exercises 38 and 39, p. 58]{Lang-algebra}, the basic problem being that
the other subgroups of $\psls$ are too small.  

We will now show that $\C{G}$ has the homotopy type of a 
three-dimensional complex and that $H_2 (\C{G})\neq 0$.
Any chain of length four (recall that a chain of length four has five vertices)
lies under a subgroup isomorphic to $S_4$, and in fact
must contain a coset $xH$ where $H\isom D_4$ or $\Z/4$.
But cosets of copies of $D_4$ and $\Z/4$ may be removed from $\C{G}$
without changing the homotopy type (apply Quillen's Theorem A
to the inclusion map).  Since the reduced Euler characteristic of $\C{G}$ is 
$17\cdot 168$~\cite[Table I]{Brown-coset-poset} and the only even dimension
in which $\C{G}$ has homology is dimension two, it now follows that $H_2 (\C{G})$ 
has rank at least $17\cdot 168$.

The method of removing cosets does not seem to show that $\C{G}$ has the 
homotopy type of a two-dimensional complex.  In fact, we expect (see  
Question~\ref{homology}) that since $H_2 (L(G))$ is non-zero,
$H_3 (\C{G})$ is non-zero as well.  (One shows that $H_2 (L(G)) \neq 0$ as follows:
We have $\tilde{\chi} (L(G)) = \mu_G (\{1\}) = 0$~\cite{Hall-eulerian}.  The above argument
shows that $L(G)$ has the homotopy type of a two-dimensional complex, and hence 
rank$\,H_1 (L(G)) =$ rank$\,H_2 (L(G))$.  
Shareshian has shown that $H_1 (L(G))\neq 0$~\cite[Lemma 3.11]{Shareshian-shell-solv}, and
in fact Shareshian's argument shows that rank$\,H_1 (L(G))\geqs 21$.)

\section{Connectivity of the Subgroup Poset}$\label{L(G)}$

In this section we will consider the connectivity of the subgroup poset, focusing in particular
on $\pi_1 (L(G))$ for finite groups $G$.  For finite solvable groups, we have the following theorem of Kratzer
and Th\'{e}venaz~\cite{Kratzer-subgroups}, which is strikingly similar to 
Theorem~\ref{solvable-coset-poset}.
As mentioned in~\cite{Thevenaz-top-homology}, this result 
follows easily by induction from Bj\"{o}rner and Walker's homotopy complementation formula
(Lemma~\ref{hom-comp}).

\begin{theorem}$\label{solvable-subgroup-poset}$
Let $G$ be a finite solvable group, and let 
$$1 = N_0\normal N_1 \normal\cdots \normal N_d = G$$
be a chief series for $G$.  Let $c_i$ denote the number of complements of 
$N_i/N_{i-1}$ in $G/N_{i-1}$.  Then $L(G)$ is homotopy equivalent to a bouquet of 
$c_1\cdot c_2\cdot \cdots \cdot c_d$ spheres of dimension $d-2$.
\end{theorem}

The number $d$ in this theorem is at least $\pi (G)$, the number of distinct primes
dividing $o(G)$ (this is immediate, since each factor $N_i/N_{i-1}$ is of prime-power order).
Note that if $c_i = 0$ for some $i$, then $L(G)$ is in fact contractible.  
Actually, more is true.

\begin{lemma}$\label{uncomplemented}$
Let $G$ be a group with a normal subgroup $N$.  If $N$ does not have
a complement in $G$, then $L(G)$ is contractible.
\end{lemma}
{\bf Proof.} This is an immediate consequence of~\cite[Theorem 3.2]{Bjorner-hom-comp}.
$\hfill \Box$

\vspace{.15in} We will say that a group $G$ is complemented if each normal subgroup $N<G$
has a complement.  The above result then reduces the study of the homotopy type of the subgroup 
poset to the case of complemented groups.

\begin{definition}
For any group $G$, let $d (G)$ denote the length of a chief series for $G$.
If $G$ does not possess a chief series, then we write $d (G) = \infty$.  
\end{definition}

In other words, $d (G)$
is the rank of the lattice of normal subgroups of $G$ (ordered by inclusion).

From Theorem~\ref{solvable-subgroup-poset} we see that if $G$ is a finite, solvable,
complemented group, then $L(G)$ is $k$-connected if and only if $d (G) \geqs k+3$.  Our first
goal in this section is to show that the ``if" portion of this statement is true
for any group (it is clearly true for any non-complemented group).

We now state the homotopy complementation formula of~\cite{Bjorner-hom-comp}, specialized
to the case of the proper part of a bounded lattice.  
Let $\overline{L}$ be a bounded lattice, and let 
$L = \overline{L} - \{\hat{0}, \hat{1}\}$ be its proper part.
We say that elements $p,q$ in a poset $P$ are \emph{complements} if $p\meet q = \hat{0}$
and $p\join q = \hat{1}$, and we denote the set of complements of $p\in P$ by $p^{\bot}$.  
We say that a subset $A\subset P$ is an \emph{antichain}
if no two distinct elements of $A$ are comparable.  

\begin{lemma}[Bj\"{o}rner-Walker]$\label{hom-comp}$
Let $L$ be the proper part of a bounded lattice.  Say there is an element
$x\in L$ such that $x^{\bot}$ is an antichain.  Then $L$ is homotopy equivalent to
$$\bigvee_{y\in x^{\bot}} {\rm Susp\,} (L_{<y} * L_{>y}).$$
\end{lemma}

In this special case, Bj\"{o}rner and Walker's proof becomes surprisingly simple.

For $G$ finite and $N$ abelian, the following result appears 
as~\cite[Corollaire 4.8]{Kratzer-subgroups}.  

\begin{proposition}$\label{wedge-decomp}$
Let $G$ be a group with normal subgroup $N$.  Then $L(G)$ is homotopy equivalent to
$\bigvee_{H\in N^{\bot}}$ Susp $(L(G/N)*L_H (N))$, where $L_H (N)\subset L(N)$ denotes the
poset of $H$-invariant subgroups.
\end{proposition}
{\bf Proof.}  Note that when $N\normal G$, the lattice-theoretic and group theoretic notions
of complement coincide.  
It is easy to check that the complements of $N$ form an antichain in $L(G)$ (when $G$ is finite,
they all have the same order).  The result now follows easily from Lemma~\ref{hom-comp},
noting that any subgroup $T$ containing a complement $H\in N^{\bot}$ has the form
$T = I\rtimes H$ where $I = T\cap N$ is $H$-invariant, and that $H\isom G/N$ implies 
$L(G)_{<H}\isom L(G/N)$.
$\hfill \Box$

\vspace{.15in}
Note that the proposition implies that if $G$ is non-simple, $\pi_1 (L(G))$ is free.

Our next result strengthens and generalizes part of~\cite[Proposition 4.2]{Kratzer-subgroups}.

By convention, we will say that a space is ($-1$)-connected if and only if it is non-empty, 
and that any space is ($-2$)-connected.  In addition, we do not consider the empty space to
be 0-connected.

\begin{proposition}$\label{subgroup-k-connectivity}$
Let $G$ be a group with normal subgroup $N\normal G$.  Then if $L(G/N)$ is $k$-connected,
$L(G)$ is $(k+1)$-connected.  In particular, if $L(G/N)$ is contractible, so is $L(G)$.
\end{proposition}
{\bf Proof.} If $L(G/N)$ is $k$-connected, then its join with any space $X$ is 
$k$-connected and Susp $(L(G/N) * X)$
is ($k+1$)-connected~\cite{Milnor-univ-bundles-2}.  
Thus Proposition~\ref{wedge-decomp} shows that $L(G)$ is ($k+1$)-connected.
$\hfill \Box$

\begin{theorem}$\label{normal-length}$
For any group $G$, $L(G)$ in ($d (G) - 3$)-connected.  In particular, if $d(G) = \infty$,
then $L(G)$ is contractible.
\end{theorem}
{\bf Proof.} When $d (G)$ is finite, this follows immediately from 
Proposition~\ref{subgroup-k-connectivity} by induction.
When $d (G) = \infty$, a similar argument shows that $L(G)$ is $k$-connected for any $k$
and hence contractible.
$\hfill \Box$

\vspace{.15in}
In order to improve on the above result, it is necessary to understand which groups have
$L(G)$ path-connected.  For non-simple finite groups, we have the following
result.

\begin{lemma}$\label{subgroup-connectedness}$
Let $G$ be a non-simple finite group.  Then $L(G)$ is disconnected if and only if 
$G\isom A\rtimes \Z/p$, where $A\neq \{1\}$ is elementary abelian, $p$ is prime, and 
$L_{\Z/p} (A)$ is empty.
\end{lemma}
{\bf Proof.}  If $G$ has the above form, then clearly $L(G)$ is disconnected ($\Z/p$ is an isolated
vertex).

In the other direction,
since $G$ is not simple we may choose a normal subgroup $N\in L(G)$ and
Proposition~\ref{wedge-decomp} gives 
$$L(G)\heq \bigvee_{H\in N^{\bot}} {\rm Susp\,\,} (L(G/N) * L_H (N)).$$
Thus in order for $L(G)$ to be disconnected, $N$ must have a complement $H$ such that
$L(G/N) = L(H) = \emptyset$ and $L_H (N) = \emptyset$.  
Thus $H\isom \Z/p$ for some prime $p$ and $H$
is a maximal subgroup of $G$.  If $G$ is a $p$-group, then $o(G) = p^2$ and we must have
$G\isom \Z/p\cross \Z/p$.  If $G$ is not a $p$-group, then the discussion after 
Theorem~\ref{M-V1} shows that $N$ is elementary abelian.
$\hfill \Box$

\begin{definition}  We will denote the collection of finite groups of the form $A\rtimes \Z/p$ 
(with $A$ elementary abelian, $p$ prime, and $L_{\Z/p} (A) = \emptyset$) by $\mF$.  The collection
of groups in $\mF$ with $A$ non-trivial will be denoted by $\mF'$.
\end{definition}

\begin{remark} The above proposition raises the interesting question of when a simple
group has a connected subgroup poset.  Maria Silvia Lucido [private communication] 
has shown that $L(G)$ is connected for any finite non-abelian simple group $G$.  Her
proof relies on the classification of the finite simple groups.
Assuming this result, we see that for finite groups $G$, $L(G)$ is disconnected if and only if 
$G\in \mF'$.
\end{remark}

The following result is immediate from Lemma~\ref{subgroup-connectedness} and 
Proposition~\ref{subgroup-k-connectivity}.

\begin{proposition}$\label{quotient}$
If a finite group $G$ has a non-simple quotient $\overline{G}$ which is not in $\mF$,
then $L(G)$ is simply connected.
\end{proposition}

Of course, one can deduce results about higher connectivity as well, and again (in light of
the above remark) the assumption of non-simplicity can presumably be removed.

Our next goal is to study the simple connectivity of $L(G)$ when $G$ can be written as a
direct product.  Assuming connectedness of $L(G)$ for simple groups, we in fact obtain
a characterization of direct products with simply connected subgroup posets.
The reader may derive the following analogue of Lemma~\ref{brown-direct-products} for the
subgroup poset of a direct product, using~\cite[Proposition 2.5]{Kratzer-subgroups} in place 
of~\cite[Proposition 1.9]{Quillen-p-subgroups}.  This result 
generalizes~\cite[Proposition 4.4]{Kratzer-subgroups}.

\begin{lemma}$\label{coprime-subgroup}$
If $H$ and $K$ are coprime groups, then 
$L(H\cross K)\heq$ Susp $(L(H) * L(K))$.
\end{lemma}

\begin{remark} The apparent incongruity between Lemma~\ref{coprime-subgroup}
and Proposition~\ref{wedge-decomp} may be explained as follows.  The extra terms in
the latter decomposition are contractible: these terms correspond to complements $K'\neq K$ of $H$,
and if $q_1:H\cross K\to H$ is the projection map, 
coprimality implies that $I\leqs I\cdot q_1 (K')\geqs q_1 (K')$ is a conical contraction of 
$L_{K'} (H)$.
\end{remark}

\begin{proposition}$\label{pi-subgroup-direct-product}$
Let $G$ be a finite group which admits a non-trivial direct-product 
decomposition.  If $L(G)$ is not simply connected, then either 
$G\isom H\cross \Z/q$ ($q$ prime) with $H\in \mF'$ or $G\isom S_1\cross S_2$ where $S_1$ and 
$S_2$ are simple groups.  In the former case, $\pi_1 (L(G))\neq 1$.
\end{proposition}
{\bf Proof.}  Let $G = H\cross K$.
Since $G$ has quotients $H$ and $K$, Proposition~\ref{quotient} shows that if $\pi_1 (L(G))\neq 1$,
then $H$ and $K$ are either simple or in $\mF$.

If $H$ and $K$ are both in $\mF'$,
then we may write $H = A_H \rtimes \Z/p$, $K = A_K\rtimes \Z/q$ ($p,q$ prime; $A_H, A_K$ 
non-trivial elementary abelian groups) and we have a chief series
$$1\normal A_H\normal H\normal H\cross A_K\normal G.$$
Theorem~\ref{normal-length} now shows that $\pi_1 (L(G)) = 1$.

Next, say that $H$ is simple and $K\in \mF'$.
If $o(H)$ is not prime, then $H$ and $K$ are coprime and Lemma~\ref{coprime-subgroup} shows that
$L(G)\heq$ Susp $(L(H)*L(K))$.  Since $L(H)$, $L(K)\neq \emptyset$, we have $\pi_1 (L(G)) = 1$.
If $H\isom \Z/p$ ($p$ prime), then the term Susp $(L_K (H) * L(K)) =$ Susp $(L(K))$ 
appears in the wedge decomposition for $L(G)$, and since $L(K)$ is disconnected $L(G)$ is not 
simply connected.
$\hfill \Box$

\begin{remark} Assuming Maria Silvia Lucido's result that $L(G)$ is connected when $G$ is a finite
non-abelian simple group,
Proposition~\ref{wedge-decomp} shows that any finite direct product $S\cross H$ 
($S$ a finite non-abelian simple group) has a simply connected subgroup poset.
In addition, $L(\Z/p\cross \Z/q)$ is clearly disconnected ($p,q$ prime), and
hence we see that if $H$ and $K$ are non-trivial finite groups, $L(H\cross K)$ is simply connected
unless $H\in \mF$ and $o(K)$ is prime
(or vice versa).
\end{remark}

We finish by examining those finite groups $G$ for which we have yet to determine 
whether or not $L(G)$ is simply connected.  From now on we assume that $L(S)$ is connected
when $S$ is a finite simple group.

Our results do not apply to non-abelian simple 
groups $G$, but Shareshian has shown that if $G$ is a 
minimal simple group (i.e. if $G$ is a finite non-abelian simple group, all of whose proper 
subgroups are solvable) then $H_1 (L(G))\neq 0$~\cite[Proposition 3.14]{Shareshian-shell-solv}.  

Any non-simple (complemented) finite group for which we have not determined simple connectivity 
of $L(G)$ may be written as (non-trivial) semi-direct
product $H\rtimes K$, where $K\in \mF$.  We break these
groups down into two cases.  First we consider groups in which $K\in \mF'$, 
and then we consider groups in which every proper, 
non-trivial quotient has prime order.  The following lemma will be useful, and is a
simple consequence of~\cite[Theorem 7.8]{Isaacs-algebra}.

\begin{lemma}$\label{simple-direct-power}$
If $N\normal G$ is a minimal normal subgroup, then $N$ is a direct power of a simple group.
In particular, if $G$ is a group in which the only characteristic subgroups are $1$ and $G$, then 
$G$ is a direct power of a simple group ($G$ is a minimal normal subgroup of $G\rtimes Aut (G)$).  

If $S$ is a finite, non-abelian simple group,
then in the direct power $S^n$, the $n$ standard copies of $S$ are permuted amongst
themselves by $Aut (S^n)$.
\end{lemma} 

Let $G = H\rtimes K$ with $K= A\rtimes \Z/p\in\mF'$.
We may assume that $G$ is not also in $\mF$.
If $H$ has a complement $K'$ such that $L_{K'} (H) = \emptyset$, then the wedge decomposition
for $L(G)$ contains the factor Susp $(L_{K'} (H)* L(K')) =$ Susp $(L(K'))$ and since
$K'\isom K$ we see that $L(K')$ is disconnected and $\pi_1 (L(G))\neq 1$ (in fact it is
not difficult to check that, if $k$ is the number of complements $K'\in H^{\bot}$ with 
$L_{K'} (H) = \emptyset$, then
$\pi_1 (L(G))$ is a free group on $k(1+ o(A))$ generators).  (Note that the groups described
in Proposition~\ref{pi-subgroup-direct-product} for which $\pi_1 (L(G))\neq 1$ are also
of this form.)
On the other hand, if $L_{K'} (H)\neq \emptyset$ for each
$K'\in H^{\bot}$, then every term in the wedge decomposition for $L(G)$ is simply connected,
so $L(G)$ is simply connected.  So in this case, $\pi_1 (L(G))\neq 1$ if and only if some 
complement of $H$ a is maximal subgroup of $G$.

Finally, we consider finite, non-simple (complemented) groups $G$ in which every proper, non-trivial quotient 
has prime order.  
Up to isomorphism, any such group may be written as $G = H\rtimes \Z/p$, and $H$ must be a 
minimal normal subgroup of 
$G$.  The lemma now shows that $H\isom S^n$ with $S$ simple.  If $o(S)$ is prime, then
$G\in \mF'$ so we assume that $S$ is a non-abelian simple group.
Now, $\Z/p\,$ acts on the $n$ standard copies of $S$ in $S^n$ and if this action has more than
one orbit then $G$ has a quotient which is not of prime order 
(and in fact $L(G)$ is simply connected by Proposition~\ref{quotient}).  So we need only consider
the cases $G = S\rtimes \Z/p$ and $G = S^p\rtimes \Z/p$, with $\Z/p$ inducing exactly 
one orbit on the standard copies of $S$.  In the latter case, $\pi_1 (L(G)) = 1$:

\begin{proposition} Let $S$ be a finite, non-abelian simple group and let $p$ be a prime.
If $G = S^p\rtimes \Z/p$, with $\Z/p$ inducing exactly 
one orbit on the standard copies of $S$, then $L(G)$ is simply connected.
\end{proposition}
{\bf Proof.}  Since $L(\Z/p)$ is empty, Proposition~\ref{wedge-decomp} shows that 
$$L(G)\heq \bigvee_{K\in (S^p)^{\bot}} {\rm Susp\,} (L_K (S^p)).$$  
Hence $\pi_1 (L(G)) = 1$ if and only
if $L_K (S^p)$ is connected for all $K\in (S^p)^{\bot}$.  Say $K\in (S^p)^{\bot}$.
Let $K = \gen{k}$ and let $\phi$ be the automorphism of $S^p$ induced by $k$.  
Also, let $S_1, \ldots, S_p$ denote the standard copies of $S$ in $S^p$ 
(so that $S^p = S_1\cross \cdots \cross S_p$).
We may assume without loss of generality that $\phi (S_i) = S_{i+1}$
($1\leqs i\leqs p-1$) and $\phi (S_p) = S_1$.  

Let $D_K\in L_K (S^p)$ denote the subgroup consisting of all elements fixed by $\phi$ 
(and hence by $K$).
For each $I\in L_K (S^p)$ we will construct a path (in $L_K (S^p)$) from $I$ to $D_K$.
Let $f_i: S^p\to S_i$ denote the $i$th projection map.  It is not hard to check that
$I\in L_K (S^p)$ implies $\hat{I} = f_1 (I)\cross \cdots \cross f_p (I)\in L_K (S^p)$
(assuming $\hat{I}\neq S^p$).  Now, since $I$ is non-trivial and $K$-invariant, 
there is a non-trivial element $s = (s_1, 1, \ldots, 1)\in f_1 (I)< \hat{I}$ and since 
$\hat{I}$ is also $K$-invariant, $\phi^i (s)\in \hat{I}$ for each $i$.
Moreover, the product $s' = \prod_{i=0}^{p-1} \phi^i (s)\in \hat{I}$ is non-trivial and invariant
under $\phi$.  Thus when $\hat{I}\neq S^p$ we have a path
$I\leqs \hat{I}\geqs \gen{s'} < D_K$.

If $\hat{I} = S^p$, then $I$ surjects onto $S_i$ ($i=1,\ldots, p$) and in particular $I$ 
is not nilpotent.  The automorphism $\phi$ either fixes $I$ pointwise or 
induces an automorphism of $I$ of prime order.  In the latter case, Thompson's theorem
on fixed-point free automorphisms implies that $\phi$ fixes some non-trivial element $i\in I$.
Hence we have a path $I\geqs \gen{i}\leqs D_K$ in $L_K (S^p)$, so $L_K (S^p)$ is connected 
and the proof is complete.
$\hfill \Box$

\vspace{.15in}
In summary, we have:

\begin{theorem} Let $G$ be a finite group which is neither simple nor a semi-direct product
$S\rtimes \Z/p$ ($S$ simple and $p$ prime), and assume further that $L(G)$ is 0-connected (i.e.
$G\notin \mF$).  Then $L(G)$ is simply connected unless $G\isom H\rtimes K$ with 
$K = A\rtimes \Z/p \in \mF'$ 
and $K$ maximal in $G$.  In this case $\pi_1 (L(G))$ is a free group on $k(1 + o(A))$ generators,
where $k$ is the number of complements of $H$ which are maximal in $G$.
\end{theorem}

\section{The Homology of $\C{G}$ and $L(G)$}

We end by discussing a relationship between the homology of the coset poset and the
homology of the subgroup poset which exists, at least, for certain groups.  This discussion
is motivated in part by Lemma~\ref{homology-surjection} and Theorem~\ref{M-V2}.

\begin{question}$\label{homology}$
If $G$ is a group, then is it true that for any $n>0$

\begin{equation}\label{hom-eq}
{\rm rank\,\,} \tilde{H}_n (L(G))\leqs {\rm rank\,\,} \tilde{H}_{n+1} (\C{G})?
\end{equation}

\end{question}

For $n=0$, the question is answered affirmatively (for any group) by Theorem~\ref{M-V2}.  
Additionally, all finite solvable groups satisfy (\ref{hom-eq}).  This follows from
Theorems~\ref{solvable-coset-poset} and~\ref{solvable-subgroup-poset}.  
We leave to the reader the easy task of checking that the number of spheres in the coset
poset is greater than the number in the subgroup poset.
In light of Lemma~\ref{uncomplemented}, any non-complemented group satisfies (\ref{hom-eq})
trivially, so (for finite groups at least) 
we may restrict our attention to non-abelian simple groups and non-trivial semi-direct products 
$H\rtimes K$ with $K$ simple.

We will now show that if $p\equiv \pm 3$ (mod 8) and $p\not\equiv \pm 1$ (mod 5)
then the simple group $PSL_2 (\F_p)$ satisfies (\ref{hom-eq}).  First, we have the following
result due to Shareshian~\cite[Lemma 3.8]{Shareshian-shell-solv}.

\begin{lemma}[Shareshian]$\label{higher-psl}$
Let $p$ be an odd prime and let $G$ be a simple group isomorphic to one of the following:
\begin{enumerate}
	\item $PSL_2 (\F_p)$ with $p\equiv \pm 3$ (mod 8) and $p\not\equiv \pm 1$ (mod 5),
	\item $PSL_2 (\F_{2^p})$,
	\item $PSL_2 (\F_{3^p})$,
	\item $Sz (2^p)$.
\end{enumerate}
Then $L(G)$ has the homotopy type of a wedge of $o(G)$ circles.
\end{lemma}

The proof of the next result is analogous to Shareshian's
proof of Lemma~\ref{higher-psl}.

\begin{lemma}$\label{c(psl)}$
If $p\equiv \pm 3$ (mod 8) and $p\not\equiv \pm 1$ (mod 5)
then $\C{PSL_2 (\F_p)}$ has the homotopy type of a two-dimensional complex.
\end{lemma}
{\bf Proof (sketch).}  Note that for $p=5$, this is just Claim~\ref{A_5-dimension}.
Let $G = PSL_2 (\F_p)$.  We begin by removing from $\C{G}$
all cosets $xH$ which are not intersections of maximal cosets, i.e. we remove all cosets
$xH$ for which $H$ is not an intersection of maximal subgroups.  The resulting poset $\mC_0$
is homotopy equivalent to $\C{G}$ by~\cite[Corollary 2.5]{Shareshian-shell-solv}.  Similarly,
the poset $L_0$ consisting of all subgroups in $L(G)$ which are intersections of maximal
subgroups is homotopy equivalent to $L(G)$.

Now,
any chain of length $k$ in the coset poset corresponds to a chain of length $k-1$ in the
subgroup poset (simply take all underlying subgroups, except for identity).  
By~\cite[Lemma 3.4]{Shareshian-shell-solv} any chain in $L_0$ has length at most two, and
hence any chain in $\mC_0$ has length at most three.
Shareshian's argument
in the proof of Lemma~\ref{higher-psl} that all two-simplices in $\ord{L_0}$ 
can be removed without changing the homotopy type also shows that all three-simplices may be
removed from $C_0$ without changing the homotopy type (when removed in the correct order, each
corresponds to an ``elementary collapse").
$\hfill \Box$

\vspace{.15in} The following computation of the Euler characteristic of $\C{PSL_2 (\F_p)}$
was provided to me by Kenneth S. Brown [private communication].

\begin{lemma}$\label{euler-psl}$
If $p\equiv \pm 3$ (mod 8) and $p\equiv \pm 2$ (mod 5), then the Euler 
characteristic of $\C{PSL_2 (\F_p)}$ is 
$o(PSL_2 (\F_p)) \left(\frac{p}{12} (p-1)(p+1)-p-4\right) + 1$.
\end{lemma}
{\bf Proof.}  For any finite group $G$, $\chi (\C{G}) = - P(G,-1) + 1$, where $P(G,s)$
is the probabilistic zeta function of $G$ (see~\cite{Brown-coset-poset}).  M\"{o}bius inversion
allows one to compute $P(G,-1)$ from the M\"{o}bius function of 
$G$~\cite[Section 2.1]{Brown-coset-poset}.  When $G = PSL_2 (\F_p)$, the M\"{o}bius function has
been calculated by Hall~\cite{Hall-eulerian}, and the reader may derive the above result.
$\hfill \Box$

\begin{proposition} If $p\equiv \pm 3$ (mod 8) and $p\not\equiv \pm 1$ (mod 5),
then $G = PSL_2 (\F_p)$ satisfies (\ref{hom-eq}).
\end{proposition}
{\bf Proof.}  By Lemmas~\ref{higher-psl} and~\ref{c(psl)}, 
it suffices to check that the rank of $H_2 (\C{G})$ is at least $o(G)$, and for $p = 5$ this 
follows from Proposition~\ref{A_5}.  We now assume $p>5$.
Since the Euler characteristic of $\C{G}$ is simply 
rank$\,H_2 (\C{G}) -$ rank$\,H_1 (\C{G})$, Lemma~\ref{euler-psl} shows that
$H_2 (\C{G})$ has rank at least 
$$o(G) \left(\frac{p}{12} (p-1)(p+1) - p - 4\right) + 1.$$  
The conditions
of the proposition force $p\geqs 11$, so $\frac{p}{12} (p-1)(p+1) - p - 4 \geqs 95$ (this is an 
increasing function of $p$).  Thus rank$\,H_2 (\C{G}) \geqs o(G)$, as desired.
$\hfill \Box$

\vspace{.15in} At least two other simple groups satisfy (\ref{hom-eq}).  These are
$PSL_2 (\F_8)$ and the Suzuki group $Sz(8)$.  The proof is analogous to that given above,
using~\cite[Table I]{Brown-coset-poset} for the computation of $P(G,-1)$ and hence $\chi (\C{G})$.
Presumably, one should be able to answer Question~\ref{homology} for all the groups listed 
in Lemma~\ref{higher-psl}.

We now show that certain direct products satisfy (\ref{hom-eq}).  In particular, given
any finite collection of non-isomorphic simple groups satisfying (\ref{hom-eq}), their 
direct product $\Pi$ also satisfies (\ref{hom-eq}), and if $G$ is a finite solvable group then 
$G\cross \Pi$ satisfies $(\ref{hom-eq})$.

\begin{proposition}$\label{coprime}$
The collection of groups satisfying (\ref{hom-eq}) is closed under coprime direct products.
\end{proposition}
{\bf Proof.} Let $H$ and $K$ be coprime groups satisfying (\ref{hom-eq}), and let $G = H\cross K$.
Recall that Lemmas~\ref{brown-direct-products} and~\ref{coprime-subgroup}
show that $\C{G}\heq \C{H}*\C{K}$ and $L(G)\heq$ Susp $(L(H)* L(K))$.
We need only consider the case in which $G$ is not solvable, and in 
light of Theorem~\ref{M-V2} we need only check condition (\ref{hom-eq}) for $n\geqs 1$.

If $L(H)$ and $L(K)$ are both empty, then $G$ is solvable and we are done.  
If $L(H)$ is empty but $L(K)$
is not, then we have $L(G) \heq$ Susp $L(K)$ and hence 
$\tilde{H}_i (L(G)) \isom \tilde{H}_{i-1} L(K)$
for $i\geqs 1$.  Letting $r_i (X)$ denote the rank of the $i$th (reduced) homology group 
of the space $X$ ($i\geqs 0$), we have (for $n\geqs 1$)
$$r_n (L(G)) = r_{n-1} (L(K)) 
             \leqs r_n (\C{K})
             \leqs \sum_{i+j = n} r_i (\C{H})\cdot r_j (\C{K})
	     = r_{n+1} (\C{G}),$$
the last equality following from~\cite[Lemma 2.1]{Milnor-univ-bundles-2}.
Of course, if $L(K)$ is empty and $L(H)$ is not, the situation is symmetric.

Now, assume that $L(H)$ and $L(K)$ are each non-empty.  Then for any $n\geqs 1$ we have
$$r_{n+1} (\C{G}) = \sum_{i+j = n} r_i (\C{H}) \cdot r_j (\C{K})$$
and
$$r_n (L(G)) = r_{n-1} (L(H)*L(K))
 = \sum_{k+l = n-2} r_k (L(H)) \cdot r_l (L(K))$$
by~\cite[Lemma 2.1]{Milnor-univ-bundles-2} (note that $L(G)$ is simply connected so
there is no problem when $n=1$).
By assumption, we have 
$r_{i-1} (L(H))\leqs r_i (\C{H})$ and $r_{i-1} (L(K))\leqs r_i (\C{K})$, 
so for each $m$ ($0\leqs m\leqs n-2$) we have 
$$r_m (L(H)) \cdot r_{n-2-m} (L(K)) \leqs r_{m+1} \C{H} \cdot r_{n-(m+1)} (\C{K})$$
and thus $r_n (L(G))\leqs r_{n+1} (\C{G})$ as desired.  
$\hfill \Box$ 

\vspace{.15in}
Th\'{e}venaz has ``found" the spheres in the subgroup
poset of a solvable group, i.e. he has provided a proof of 
Theorem~\ref{solvable-subgroup-poset} by analyzing a certain collection of spherical
subposets of $L(G)$~\cite{Thevenaz-top-homology}.  
It would be interesting to explore similar ideas in the
coset poset, and in particular such a proof of Theorem~\ref{solvable-coset-poset}
might allow one to explicitly construct an injection from 
$H_n (L(G))$ into $H_{n+1} (\C{G})$ (for $G$ finite and solvable), and could shed
further light on Question~\ref{homology}.  
\bibliography{references}

\bibliographystyle{hamsplain}

\end{document}